\documentclass[11pt]{article}

\usepackage{color}
\usepackage{latexsym}
\usepackage{amssymb}
\usepackage{amsthm}
\usepackage{amsmath}
\usepackage{graphicx}
\usepackage{enumerate}
\usepackage{hyperref}
\usepackage{color}
\usepackage{amsfonts}

\newtheorem{Theorem}{Theorem}[part]
\newtheorem{Definition}{Definition}[part]
\newtheorem{Proposition}{Proposition}[part]

\newtheorem{Corollary}{Corollary}[part]
\newtheorem{Remark}{Remark}[part]

\def \Frac{\displaystyle\frac}

\def\Esp#1{\mathbb{E}\left[#1\right]}

\def \N{\mathbb{N}}
\def \R{\mathbb{R}}

\def \Gc{{\cal G}}

\def \Ic{{\cal I}}

\def \Lc{{\cal L}}
\def \Pc{{\cal P}}

\def \Oc{{\cal O}}
 
\def \Sc{{\cal S}}

\def \Lb{{\bf L}}

\def \Eb{{\bf E}}
\def \Pb{{\bf P}}
\def \Ab{{\bf A}}

\def \eps{\varepsilon}

\def \ep{\hbox{ }\hfill$\Box$}

\def\Dt#1{\Frac{\partial #1}{\partial t}}

\def\reff#1{{\rm(\ref{#1})}}

\def\beqs{\begin{eqnarray*}}
\def\enqs{\end{eqnarray*}}
\def\beq{\begin{eqnarray}}
\def\enq{\end{eqnarray}}

\begin{document}
\title{Probabilistic Representation and Approximation \\ for Coupled Systems of Variational Inequalities\footnote{\textbf{Acknowledgement}. We would like to thank both anonymous  referees for useful comments. }
\vspace{0mm}}
\author{
\begin{tabular}{ccc}
Romuald ELIE & ~~ & Idris KHARROUBI\\
CEREMADE, CNRS, UMR 7534, & ~~ &  LPMA, CNRS, UMR 7599,  \\
Universit{\'e} Paris-Dauphine, & ~~  & Universit{\'e} Paris 7,\\
and CREST & ~~ & and CREST \\
\sf elie@ensae.fr &  & \sf kharroubi@ensae.fr
\end{tabular}
}

\date{  }

\maketitle
\begin{center}
\vspace{-9mm}
\large{May 2010 \\
First version: {September 2009}}\\[15mm]
\end{center}

\vspace{-8mm}
\begin{abstract}
Our study is dedicated to the probabilistic representation and numerical approximation of solutions of coupled systems of variational inequalities. 
We interpret the unique viscosity solution of a coupled system of variational inequalities as the solution of a one-dimensional constrained BSDE with jumps. 
This new representation allows for the introduction of a natural probabilistic numerical scheme for the resolution of these systems. 

\end{abstract}

\vspace{0mm}

\vspace{0mm}

\noindent {\bf MSC Classification (2000):}  93E20, 60H10, 60H30, 35K85, 49L25.




\section{Introduction}
\setcounter{equation}{0} \setcounter{Assumption}{0}
\setcounter{Theorem}{0} \setcounter{Proposition}{0}
\setcounter{Corollary}{0} \setcounter{Lemma}{0}
\setcounter{Definition}{0} \setcounter{Remark}{0}

Pardoux and Peng (1992) 
developed the theory of backward stochastic differential equations, providing a probabilistic representation of solutions of quasi-linear parabolic PDEs.
 Coupling the diffusion process with a pure jump process, Pardoux et al. (1997) 
 extend this representation to systems of coupled semilinear PDEs with different linear differential operators on each line. Introducing restrictions on the domain of the backward process, El Karoui et al. (1997) 
 cover the class of variational inequalities. Constraining instead the jump part of the solution, Kharroubi et al. (2010) 
 consider quasilinear variational inequalities.

The focus of this note is to extend this type of Feynman-Kac representation to the more general class of coupled systems of quasilinear variational inequalities arising, for example, in optimal impulse or switching problems. We will typically consider systems of PDE of the form 
 \beq \label{SVI}
 \Big[ - \Dt{v_i} - \Lc^i v_i - f(i,.,(v_k)_{1\le k\le m},\sigma(i,.)^\top D_xv_i)\Big]
 \wedge \min_{1\leq j\leq m} h_{}(i,j,.,v_i,v_j,\sigma(i,.)^\top D_xv_i) = 0, \\
 \mbox{ on } \Ic\times[0,T)\times \R^d, \qquad  \mbox{ with terminal condition }\quad v_i(T,.)~=~g(i,.) \mbox{ on } \Ic\times\R^d\,,\label{SVIT}
 \enq
 where, for any $i\in\Ic:=\{1,\ldots,m\}$,  $\Lc^i$ is a linear second order local operator 
  \beq \label{def Li}
  \Lc^i v_i(t,x) &:=& b(i,x) \cdot D_x v_i(t,x) + \frac{1}{2}{\rm tr}(\sigma\sigma^\top(i,x) D_x^2v_i(t,x))\,,  
  \enq 
  and $b$, $\sigma$, $f$, $h$ and $g$ are Lipschitz continuous functions. 
  As observed by Bouchard (2009), 
  this PDE appears in the resolution of optimal switching problems as well as stochastic target problems with jumps. 
 The major difficulty arises from the coupling between all the components $(v_i)_{i\leq m}$ of the solution and the use of different linear operators at each line. When $m$ is large, the numerical resolution of \reff{SVI}-\reff{SVIT} by classical PDE approximation methods is very tricky and highly computational. We intend to provide here a probabilistic representation to \reff{SVI}-\reff{SVIT} leading to an efficient probabilistic numerical scheme. When $b$ and $\sigma$ are independent of the regime $i\in\Ic$ and the constraint functions are of the form $h_{}:(i,j,.,y_i,y_j,.)\,\mapsto\, y_i-y_j-c_{i,j} $,  Hu and Tang (2007) 
 interpret the vector solution to \reff{SVI}-\reff{SVIT} as a multi-dimensional BSDE with terminal condition and oblique reflections. The challenging derivation of a convergent numerical approximation for this type of BSDE is of great interest and is currently under study. The approach of this paper relies instead on a recent reinterpretation of obliquely multi-dimensional reflected BSDEs in terms of one-dimensional constrained BSDEs with jumps, as introduced in Elie and Kharroubi (2009). 
  The idea is to consider, as in Pardoux et al. (1997), 
  a random regime driven by a pure jump transmutation process, allowing to retrieve simultaneously  some information concerning all the components of the solution.\vspace{2mm} 

Given a $d$-dimensional Brownian motion $W$ and an independent Poisson measure $\mu$ on $\R_+\times\Ic$, we consider, for any initial condition $e:=(t,i,x)\in [0,T]\times\Ic\times\R^d=:E$, the unique $\Ic\times\R^d$-valued solution $(I^e_s,X^e_s)$ of the SDE: 
  \beq\label{defXswitch}
 \left\{
 \begin{array}{lll}
 I_{s} & = & i+\int_{t}^s\int_{\Ic}(j-I_{r-})\mu(dr,dj)\\
 X_{s} & = & x+\int_{t}^sb(I_{r},X_{r})dr+\int_{t}^s\sigma(I_{r},X_{r}) \cdot dW_{r}
 \end{array}
 \right.,\quad\qquad t\le s \le T\;.
 \enq
 Formally, given a smooth solution $(v_{i})_{i\in\Ic}$ to \reff{SVI}-\reff{SVIT}, the process $Y$ $:=$ $v_{I^e}(.,X^e_.)$ satisfies
 \begin{align}
 Y_t = g_{}(I^e_T,X^e_T) + \hspace{-1mm}\int_t^T \hspace{-2mm}f_{}(I^e_s,X^e_s,Y_s + U_s, Z_s) ds  + K_T - K_t   -   \hspace{-1mm}\int_t^T \hspace{-2mm}  Z_s \cdot dW_s  -  \hspace{-1mm} \int_t^T \hspace{-2mm}\int_\Ic \hspace{-1mm}  U_s(j)   \mu(ds,dj)
 \label{BSDEgen}
 \end{align}
 on $[0,T]$, where we denote $Z_s:=\sigma_{}^\top(I^e_{s-},X^e_{s})D_xv_{I^e_{s-}}(s,X^e_{s})$, $U_s(.):= v_{.}(s,X^e_{s})- v_{I^e_{s-}}(s,X^e_{s})$,
and\\ $K_s := \int_0^s [-\Dt{v_{I^e_{u}}}$ $-\Lc^{I^e_{u}} v_{I^e_{u}} -   f_{}(I^e_{u},.,(v_k)_{1\le k\le m},\sigma^\top(I^e_{u},.) D_xv_{I^e_{u}})] (u,X^e_u) du$. 
Since $v$ satisfies \reff{SVI}, we expect the following constraint to be satisfied:
 \beq
 \label{cons}
  h_{}(I^e_{s-},j,X^e_s,Y_{s-}, Y_{s-} + U_s(j), Z_s) & \geq & 0\;,  \qquad j\in\Ic,\;\; t\le s\le T\,.
 \enq 
 The BSDE \reff{BSDEgen} combined with constraint \reff{cons} falls into the class of constrained BSDEs with jumps and admits a unique minimal solution under mild conditions on the coefficients.  We reinterpret the $Y$-component of the solution as the unique viscosity solution to the coupled system of variational inequalities \reff{SVI}-\reff{SVIT}.
This new Feynman-Kac representation is meaningful to the BSDE literature since:
 \begin{itemize}
\setlength\itemsep{0pt}
\setlength\parskip{0pt}
\setlength\leftmargin{\parindent}
\setlength\itemindent{-\parindent}
 \item It extends the results of Kharroubi et al. (2010) 
 to more general constraints and driver functions depending on $U$. This allows for a strong coupling between the dynamics of the value function components and gives a minimality condition in some particular cases.
 \item It generelizes the conclusions of Peng and Xu (2007) 
 derived in the no-jump case.
 \item It offers a PDE representation to reflected BSDEs with interconnected obstacles introduced in Hamad\`ene and Zhang (2008) 
 since they relate directly to constrained BDSE with jumps, see Elie and Kharroubi (2009). 
 \item It generalizes the use of diffusion-transmutation process in Pardoux et al. (1997) 
  to systems of variational inequalities.
 \end{itemize}
 This representation leads to a natural probabilistic algorithm for the resolution of \reff{SVI}-\reff{SVIT}.
 The constrained BSDE with jumps is replaced by a penalized BSDE with jumps, which is approximated by the discrete-time scheme studied in Bouchard and Elie (2008) 
 and Gobet et al. (2006).
 This leads to a convergent numerical scheme based on time discretization, Monte Carlo simulations and projections.\vspace{1mm}
 
 The rest of the paper is organized as follows. In Section 2, we discuss existence, uniqueness and penalization, and give a minimality condition for constrained BSDEs with jumps \reff{BSDEgen}-\reff{cons}. Section 3 presents the viscosity properties and the numerical approximation is detailed in the last section.

\vspace{-2mm}

\paragraph{Notation.} 
Throughout this paper, we are given a finite horizon $T$ and a probability space $(\Omega,\Gc,\Pb)$ endowed with a $d$-dimensional standard Brownian motion $W$ $=$ $(W_t)_{t\geq 0}$, and an independent Poisson random measure $\mu$ on $\R_+\times \Ic$,  with intensity measure $\lambda(di)dt$ for some positive finite measure $\lambda$ on $\Ic$ $:=$ $\{1,\ldots,m\}$. We denote $E:=[0,T]\times\Ic\times\R^d$.
For a smooth function $\varphi:~[0,T]\times\mathbb{R}^d\times\Ic\rightarrow\mathbb{R}$, $\frac{\partial \varphi}{\partial t}$, $D_x \varphi$ and $D^2_x\varphi$ denote resp. the derivative of $\varphi$ w.r.t. $t$, the gradient and the Hessian matrix of $\varphi$ w.r.t. $x$. The dependence in $\omega\in\Omega$ is omitted when it is obvious.


\section{Constrained Forward Backward SDEs with jumps}\label{SecBSDE}
\setcounter{equation}{0} \setcounter{Assumption}{0}
\setcounter{Theorem}{0} \setcounter{Proposition}{0}
\setcounter{Corollary}{0} \setcounter{Lemma}{0}
\setcounter{Definition}{0} \setcounter{Remark}{0}

We present in this section the constrained forward backward SDEs with jumps and  recall the existence and uniqueness results of Elie and Kharroubi (2009). 
We 
discuss the correspondence between the value function associated to $Y$ and the $U$ component of the solution. Under additional regularity of the value function, we provide a Skorohod type minimality condition for the considered BSDE.

\vspace{-2mm}

\subsection{Existence and uniqueness of a minimal solution via penalization}\label{SubSecPenalization}


 As discussed above, the forward process is a transmutation-diffusion process composed of a pure jump process $I$ and a diffusion without jump $X$ whose dynamics depends on $I$. For any initial condition $e:=(t,i,x)\in E$, $(I^e,X^e)$ is the unique solution to \reff{defXswitch} starting from $(i,x)$ at time $t$. \\

 For any initial condition $e\in E$, a solution to the constrained BSDE with jumps is a quadruplet $(Y^e,Z^e,U^e,K^e)\in\Sc^2\times\Lb_W^2\times\Lb^2_{\tilde\mu}\times\Ab^2$ satisfying \reff{BSDEgen}-\reff{cons}, where 
  \begin{itemize}
\setlength\itemsep{0pt}
\setlength\leftmargin{\parindent}
\setlength\itemindent{-\parindent}
\item $\Sc^2$ is the set of real valued $\Gc$-adapted  {\sl c\`adl\`ag} processes $Y$ on $[0,T]$ s.t.\\
    $
    \|Y\|_{\Sc^2} := \Esp{\sup_{0\le r \le T} |Y_r|^2}^{\frac{1}{2}} <\infty
    $,
 \item $\Lb^p_W$ is the set of progressive $\R^d$-valued processes $Z$ s.t.
    $
    \|Z\|_{\Lb^p_W} :=  \Esp{\left(\int_0^T |Z_r|^p dr\right)}^{\frac{1}{p}} <\infty
    $, $p$ $\geq$ $1$,
 \item $\Lb^p_{\tilde\mu}$ is the set of $\Pc\otimes\sigma(\Ic)$ measurable maps $U$ $:$ $\Omega\times [0,T]\times \Ic \rightarrow \R$ s.t.\\
    $
    \|U\|_{\Lb^2_{\tilde\mu}}:= \Esp{ \int_0^T \int_\Ic |U_s(j)|^2 \lambda(dj) ds }^{\frac{1}{p}}
    <\infty
    $, $p$ $\geq$ $1$,
 \item ${\Ab^2}$ is the  closed subset of ${\Sc^2}$ composed by nondecreasing processes $K$ with $K_0$ $=$  $0$.
\end{itemize}
Furthermore, $(Y,Z,U,K)$ is referred to as the minimal solution to \reff{BSDEgen}-\reff{cons} whenever we have $Y\leq Y'$ a.s., for any other solution $(Y',Z',U',K')$
. In order to ensure existence and uniqueness of a minimal solution to \reff{BSDEgen}-\reff{cons} for any initial condition, we make the following assumptions.\vspace{2mm}

\noindent \textbf{(H0)} The following holds:\vspace{-2mm}
 \begin{enumerate}[(i)]
\setlength\itemsep{0pt}
\setlength\leftmargin{\parindent}
\setlength\itemindent{-\parindent}
\item 
There exists a constant $L$ s.t.
\beqs 
 |f(i,x,(u_{j})_{j\in\Ic},z)-f(i,x,(u_{j}')_{j\in\Ic},z')|  & \leq &  L |(z,(u_{j})_{j\in\Ic})-(z',(u_{j}')_{j\in\Ic})|\;,
 \nonumber\\
  |h_{}(i,j,x,y,u_j,z,j)-h(i,j,x,y',u'_j,z')| & \leq &  L |(y,z,u_j)-(y,z',u'_j)|\,,
  \nonumber
\enqs
 for all  $(x,i,j,y,z,u,y',z',u')\in\R^d\times\Ic^2\times[\R\times\R^{d}\times\R^{\Ic}]^2$, and 
 \beqs
 |f(i,x,(u_{j})_{j\in\Ic},z)|+  |h_{}(i,j,x,y,u_j,z)| & \leq  & L\big(1+|(y,z,(u_{j})_{j\in\Ic})|\big)\;,
 \enqs
  for all  $(x,i,j,y,z,(u_{i})_{i\in\Ic})\in\R^d\times\Ic^2\times\R\times\R^{d}\times\R^{\Ic}$.
\item The function $h_{}(i,j,x,y,.,z)$ is non-increasing for all $(i,x,y,z,j)$ $\in$ $\Ic\times\R^d\times\R\times\R^d\times\Ic$.
 \item There exist two constants
$C_{1}\geq C_{2}>-1$ and a measurable map $\gamma:
\Ic\times\R^d\times\R\times\R^d\times[\R^\Ic]^2\times\Ic\rightarrow[C_2,C_1]$ such that, for any $(i,x,y,z,u,u')$ $\in$ $\Ic\times\R^d\times\R\times\R^d\times[\R^\Ic]^2$,
 \beqs
 f(i,x,y+u,z)-f(i,x,y+u',z')
 &\leq&
 \int_{\Ic}(u_{j}-u'_{j})\gamma(i,x,y,z,u,u',j)\lambda(dj)\,.
 \enqs
\end{enumerate}


\noindent {\bf (H1)} 
 For any $e=(t,i,x)\in E$, there exists a quadruple $(\tilde Y^e,\tilde Z^e,\tilde U^e,\tilde K^e)$ $\in$ ${\bf \Sc^2}\times{\Lb_W^2}\times{\Lb_{\tilde\mu}^2}\times{\Ab^2}$ solution to \reff{BSDEgen}-\reff{cons}, with $\tilde Y^e_t$ $=$ $\tilde v_{I^e_t}(t,X^e_t)$, for some deterministic function $\tilde v$ satisfying 
$|\tilde v_i(t,x)|  \leq  C(1+|x|)$ on $E$.


\vspace{1mm}

\noindent We provide in Remark \ref{RemSufficientH1} a more tractable sufficient condition under which {\bf (H1)} holds. \\[-2mm]

The construction of the minimal solution is done by penalization. 
 For any initial condition $e\in E$ and $n\in \N$, we introduce $(Y^{e,n},Z^{e,n},U^{e,n})$ solution to the following penalized BSDE 
 \beq\label{BSDEpen}
 Y_t & = & g_{}(I^e_{T},X^e_{T}) + \int_t^T f_{}(I^e_{s},X^e_{s},Y_s+U_s,Z_s) ds - \int_t^T \int_\Ic   U_s(j)  \mu(ds,dj) - \int_t^T Z_s \cdot dW_s \nonumber\\
 &+&  n \int_t^T  \int_\Ic [h_{}(I^e_{s-},j,X^e_{s},Y_{s-},Y_{s-}+U_s(j),Z_s)]^-\lambda(dj)ds\,,\qquad 0\le t \le T\,.
 \enq
 Under {\bf (H0)}, we get from Barles et al. (1997) 
 existence and uniqueness of a solution of \reff{BSDEpen}. We introduce $K^{e,n}:=\int_0^.\int_\Ic [h_{}(I^e_{s-},j,X^e_{s},Y^{e,n}_{s-},Y^{e,n}_{s-}+U^{e,n}_s(j),Z^{e,n}_s)]^-\lambda(dj)ds$, for any $(e,n)\in E\times \N$. 
 \begin{Theorem}\label{ThmExistsminimal}
 Suppose {\bf (H0)}-{\bf (H1)} holds.
 For any $e:=(t,i,x)\in E$, there exists a unique quadruple $(Y^e,Z^e,U^e,K^e)$ $\in$ ${\bf \Sc^2}\times{\Lb_W^2}\times{\Lb_{\tilde\mu}^2}\times{\Ab^2}$
 minimal solution to \reff{BSDEgen}-\reff{cons} with $K^e$ predictable, and $v_i:(t,x)\mapsto Y^{t,i,x}_{t}$ defines a deterministic map from $E$ into $\R$. Moreover $(Y^{e},Z^{e},U^{e})$ is the limit of the $(Y^{e,n},Z^{e,n},U^{e,n})_{n\in\N}$  in the following sense 
 \beqs 
 \|Y^{e,n}-Y^e\|_{_{{\Lb^2_{W}}}} + 
\|Z^{e,n}-Z^e\|_{_{{ \Lb^p_{W}}}} + \|U^{e,n}-U^e\|_{_{{\Lb^p_{\tilde\mu}}}} 
   & \longrightarrow & 0, ~~~n\rightarrow\infty, \qquad  1\le p < 2 \;.
 \enqs 
 \end{Theorem}

\textbf{Proof.} This result is a direct application of Theorem 2.1 in Elie and Kharroubi (2009). \ep
\vspace{2mm}
Under additional regularity on $Y^e$, we can improve the previous convergence up to $p$ $=$ $2$.
 \begin{Proposition} \label{PropConvPenError}
 If {\bf (H0)}-{\bf (H1)} holds, $(Y^{e,n})_{n\in\N}$ converges increasingly to $Y^e$, for any $e\in E$. Additionally, if the process $Y^{e}$ is quasi-left continuous in time, we have
 \beq\label{ConvPenError2}
  \|Y^e-Y^{e,n}\|_{_{\Sc^2}}   +  \|Z^e-Z^{e,n}\|_{_{\Lb^2_W}}  +  \|U^e-U^{e,n}\|_{_{\Lb^2_{\tilde\mu}}} 
  +  \|K^e-K^{e,n}\|_{_{\Sc^2}}
  \mathop{{\longrightarrow}}\limits_{n\rightarrow\infty}   0\,, \;\;\;e\in E \;.\quad
 \enq
 \end{Proposition}

 \textbf{Proof.}
 Fix $e\in E$ and observe from Proposition 2.1 in Elie and Kharroubi (2009) 
  that $Y^{e,n}$ converges increasingly to $Y^e$. Since $\mu$ is a Poisson measure, 
the process  $Y^{e,n}$ is quasi-left continuous. If $Y^e$ has the same regularity, the predictable projections of $Y^{e}$ and $Y^{e,n}$ are simply given by $(Y_{t-}^{e})_t$ and $(Y_{t-}^{e,n})_t$. This leads to $Y^{e}_{t-}=\lim_{n\to\infty} Y_{t-}^{e,n}$. We deduce from the weak version of Dini's theorem, see Dellacherie and Meyer (1980) 
 p. 202, that $Y^{e,n}$ converges uniformly to $Y^{e}$ on $[0,T]$, and the dominated convergence theorem gives us $ \|Y^e-Y^{e,n}\|_{_{\Sc^2}}  \mathop{{\longrightarrow}}\limits_{n\rightarrow\infty}  0$. 
 Combined with standard estimates of the form
\beqs
 \|Z^{e,n+p}-Z^{e,n}\|_{_{{\Lb_W^2}}}^2 + \|U^{e,n+p}-U^{e,n}\|_{_{{\Lb_{\tilde\mu}^2}}}^2
 +\|K^{e,n+p}-K^{e,n}\|_{_{{\Sc^2}}}^2 
 &\leq & C\|Y^{e,n+p}-Y^{e,n}\|_{_{{\Sc^2}}}^2,
 \enqs
 this implies that the sequences $(Z^n)$, ${(U^n)}$ and $(K^{n})$ are Cauchy and hence convergent. 
 \ep


 \begin{Remark}\label{RemConvPenError}{\rm 
 Under the additional Assumption {\bf (H2)} below, $(v_i)_{i\in\Ic}$ is interpreted as the unique viscosity solution to \reff{SVI}-\reff{SVIT}, see Theorem \ref{theouni}. In this case, $(v_i)_{i\in\Ic}$ is continuous, $Y_t=v_{I_t}(t,X_t)$ is quasi-left continuous and Proposition \ref{PropConvPenError} holds.}
 \end{Remark}

 We denote by $(v^n)_{n\in\N}$ the sequence of deterministic functions defined by $v^n \,:\; e\in E \mapsto Y_t^{e,n}$ and we shall use indifferently the notation $v^n(t,i,x)$ or $v^n_i(t,x)$, for $(t,i,x)\in E$. 
 Under {\bf (H0)-(H1)}, we know from Proposition \ref{PropConvPenError} that $v$ is the pointwise limit of $(v^{n})_{n\in \N}$.



\subsection{Representation of $U$ and the minimality condition}\label{SubSeclYUX}

\begin{Proposition}\label{PropRepresU}
 Let {\bf (H0)-(H1)} holds. 
 For any $e\in E$ and stopping time $\theta$ valued in $[t,T]$, we have $Y^e_{\theta}  =  v_{I^e_\theta}(\theta,X^e_{\theta})$, and the process $U$ represents as
\beq\label{idvU}
 U_{s}^{e}(j)  & = & v_{j}(s,X^e_s)-v_{I^e_{s-}}(s,X^e_s)\,, \qquad j\in\Ic, \;\;\;\;\;  t\le s \le T \;. 
\enq
\end{Proposition}
{\bf Proof.} 
According to Proposition \ref{PropConvPenError}, we simply need to provide similar representations for the penalized BSDE \reff{BSDEpen}. 
Fix $e\in E$. For any stopping time $\theta$ valued in $[t,T]$, uniqueness of solution of \reff{BSDEpen} and the Markov property of $(I^e,X^e)$ directly give to $Y^{e,n}_{\theta}  =  v^n_{I^e_\theta}(\theta,X^e_{\theta})$. Denoting $\tilde{U}_s^{e,n}(j) := v^n_j(s,X^e_s) -v^n_{I^e_{s-}}(s,X^e_s)$, for $j\in\Ic$ and $0\le s\le T$, we deduce from \reff{BSDEpen} that
\vspace{-1mm}
\beqs
\int_\Ic \tilde U^{e,n}_s(j)\mu(ds,dj) & = &  Y^{e,n}_s-Y^{e,n}_{s-} ~ = ~ \int_\Ic \tilde U^{e,n}_s(j)\mu(ds,dj)\;, \qquad 0\le s\le T\,. 
\vspace{-1mm}
 \enqs
Therefore $\mathbb{E}\left[\int_0^T\int_\Ic(U^{e,n}_s(j)-U^{e,n}_s(j))^2\lambda(dj)ds\right]=0$ and the proof is complete.\ep

\vspace{2mm}

Under an extra regularity assumption on the function $v$ satisfied under Assumption {\bf (H2)} below, the previous representation  
 leads to a Skorohod type minimality condition for \reff{BSDEgen}-\reff{cons}.
 
 \begin{Corollary}\label{CorMinimality}
Let {\bf (H0)-(H1)} holds. Suppose $(v_i)_{i\in\Ic}$ is continuous and the function $h$ does not depend on $z$. Then, for any $e\in E$, the minimal solution $(Y^{e},Z^{e},U^{e},K^{e})$ satisfies 
\vspace{-2mm}
\beq\label{SkC}
\int_{t}^T\min_{j\in\Ic}\Big[h_{}(I^{e}_{s-},j,X^e_{s},Y^{e}_{s-},Y^{e}_{s-}+U^{e}_{s}(j))\Big]dK^{e}_{s} & = & 0.
\enq
\end{Corollary}
\textbf{Proof.}
 Fix $e\in E$. Since $(v_i)_{i\in\Ic}$ is continuous, the process $Y^{e}$ inherits the quasi-left continuity of $(I^{e},X^{e})$. Combining \reff{idvU} and Proposition \ref{PropConvPenError} leads to
$ \max_{j\in\Ic} \|U^{e}(j)-U^{e,n}(j)\|_{_{\Sc^2}} \mathop{{\longrightarrow}}\limits_{n\rightarrow\infty}  0$. 
We deduce from \reff{ConvPenError2} and Lemma 5.8 in Gegoux-Petit and Pardoux (1995), 
which also holds for c\`agl\`ad functions, that
\vspace{-2mm}
\beqs
\int_{t}^T\hspace{-1mm}\min_{j\in\Ic}\Big[h_{}(I^e_{s-},j,X^e_{s},Y^{e,n}_{s-},Y^{e,n}_{s-}+U^{e,n}_{s}(j))\Big]dK^{e,n}_{s} 
 \mathop{{\longrightarrow}}\limits_{n\rightarrow\infty} 
\int_{t}^T\hspace{-1mm}\min_{j\in\Ic}\Big[h_{}(I^{e}_{s-},j,X^e_{s},Y^{e}_{s-},Y^{e}_{s-}+U^{e}_{s}(j))\Big]dK^{e}_{s}\;.
\enqs
\vspace{-1mm}
Since 
$
\int_{t}^T\min_{j\in\Ic}\Big[h_{}(I^e_{s-},j,X^e_{s},Y^{e,n}_{s-},Y^{e,n}_{s-}+U^{e,n}_{s}(j))\Big]dK^{e,n}_{s} \leq 0
$ and  \reff{cons} holds, we get \reff{SkC}.
 \ep


\section{Link with coupled systems of variational inequalities}
\setcounter{equation}{0} \setcounter{Assumption}{0}
\setcounter{Theorem}{0} \setcounter{Proposition}{0}
\setcounter{Corollary}{0} \setcounter{Lemma}{0}
\setcounter{Definition}{0} \setcounter{Remark}{0}

In this section, we interpret the minimal solution of \reff{BSDEgen}-\reff{cons} as the unique viscosity solution of the PDE \reff{SVI}-\reff{SVIT}, thus generalizing the representation derived in Kharroubi et al. (2010), 
Pardoux et al. (1997) and Peng and Xu (2007). 

\subsection{Viscosity properties of the penalized BSDE}

The penalized parabolic integral partial differential equation (IPDE) associated to  \reff{BSDEpen} is naturally defined for each $n$ $\in$ $\N$ by
 \beq \label{IPDE}
 \left\{
 \begin{array}{ll}
 & - \Dt{\varphi_i} - \Lc^i \varphi_i  - f(i,.,(\varphi_j)_{j\in\Ic},\sigma^\top(i,.) D_x \varphi_i) 
 - n \int_{\Ic} \big[{h}_{}\big(i,j,.,\varphi_i,\varphi_j, \sigma_i^\top D_x \varphi_i \big)\big]^- \lambda(dj)   =   0 \\ 
 & \mbox{ on } [0,T)\times\R^d\times\Ic, \qquad \mbox{ and } \quad  v_i(T,.) = g(i,.) \mbox{ on } \Ic\times\R^d\,,
 \end{array}
 \right.
 \enq
 where $\Lc$ is the $m$-dimensional Dynkin operator associated to $X$, defined in \reff{def Li}.
 Since the penalized BSDE falls into the class of BSDE with jumps  studied by Pardoux et al. (1997), 
 we deduce the following Feynman-Kac representation result.

\begin{Proposition} \label{provn}
Under {\bf (H0)-(H1)}, the functions $(v^n)_n$ are continuous viscosity solutions of \reff{IPDE}. Indeed, for any $n\in\N$, $v^n(T,.)=g$ and, for any $(i,t,x)$ $\in$ $\Ic\times[0,T)\times\R^d$ and $\varphi\in C^{1,2}([0,T]\times\mathbb{R}^d)$ such that $(t,x)$ is a global minimum (resp. max.imum)  of $(v_i^n-\varphi)$, we have
 \beqs
\left[ \hspace{-0.5mm}- \Dt{\varphi} -  \Lc^i\varphi  - \hspace{-0.5mm} {f}(i,.,(v^n_j)_{j\in\Ic},\sigma^\top(i,.) D_x\varphi) - n \hspace{-1mm}\int_{\Ic}\hspace{-0.5mm}  [h_{}(i,j,.,v^n_i,v^n_j,\sigma^\top(i,.) D_x \varphi)]^- \lambda(dj)\right] (t,x) \;\geq\;   \,(\mbox{resp. }   \leq\,)  \;0.
 \enqs
\end{Proposition}

\textbf{Proof.} Fix $n\in\N$. The continuity of $v^n$ follows from similar arguments as in the proof of Lemma 2.1 in Pardoux et al. (1997). 
According to the representations detailed in the proof of Proposition \ref{PropRepresU}, the viscosity property of $v^n$ fits in the framework of Theorem 4.1 in Pardoux et al. (1997), 
up to the comparison theorem for BSDE, which is replaced by Theorem 2.5 in Royer (2006).
\ep


\subsection{Viscosity properties of the constrained BSDE with jumps}

 Formally, passing to the limit in \reff{IPDE} when $n$ goes to infinity, we expect $v$ to be a solution of \reff{SVI} on $[0,T)\times\R^d\times\Ic$. 
 As for the boundary condition, we cannot expect to have $v(T^- , .) = g$, and we shall consider the relaxed boundary condition given by 
 \beq \label{SVITT} 
 \min \Big[ v_i - g(i,.) \; , \;  \min_{j\in \Ic} h_{}\left(i,j,., v_i, v_j,\sigma^\top(i,.) D_x v_i\right)  \Big](T^- ,x) &=& 0 \mbox{ on } \Ic\times\R^d.\;\;\;
 \enq

 \begin{Remark}\label{RemSwitch}{\rm
 In the particular case where the driver function $f$ is independent of $(y,z,u)$ and the constraint function
 is given by $\tilde h_{}:(i,j,x,y,y+v,z)\mapsto -c_{i,j}-v$ with $c$ a given cost function, we retrieve the system of variational inequalities associated to switching problems
 \beq
 &\min\Big[ - \Dt{v_i} - \Lc^i v_i - f(i,.), \; \min_{j\in \Ic} \left[ v_i -v_j - c_{i,j}\right]   \Big] &= 0\,,\mbox{ on } [0,T)\times\R^d\times\Ic\;, \label{SVISwitch}\\
 &\min\Big[ v_i - g(i,.)  , \; \min_{j\in \Ic} \big[ v_i-  v_j-c_{i,j}\big] \Big](T^- ,.) & =  0\,,   \mbox{ on } \R^d\times\Ic\;. \label{SVISwitchT}
 \enq
 Thus, if \reff{SVISwitchT} satisfies a comparison theorem, $v(T^-,.)$ is the smallest function greater than $g$ satisfying \reff{SVISwitchT}. 
 In particular, we retrieve the terminal condition $v(T^-,.)=g$ proposed by Hu and Tang (2007) 
 when the terminal condition $g$ satisfies the cost constraint.
 }\end{Remark}


In order to define viscosity solutions of \reff{SVI}-\reff{SVITT}, we introduce, for any locally bounded vector function $(u_i)_{i\in\Ic}$ on
 $[0,T]\times\R^d$ its lower semicontinuous and upper semicontinuous (lsc and usc for short) envelopes $u_*$ and $u^*$
defined for $(t,x)\in[0,T]\times\R^d$ by
 \beqs
 u_*(t,x)  =  \liminf_{\tiny{(t',x')\rightarrow (t,x), t'<T}} u(t',x'), &\mbox{ and }&  u^*(t,x) =  \limsup_{\tiny{(t',x')\rightarrow (t,x), t'<T}} u(t',x')\,.
 \enqs

\begin{Definition}

\noindent A vector function $(u_i)_{i\in\Ic}$, lsc (resp. usc) on $[0,T)\times\R^d$, is called a viscosity supersolution (resp. subsolution) to \reff{SVI}-\reff{SVITT} if, for each $(i,t,x)$ $\in$ $\Ic\times[0,T]\times\R^d$ and $\varphi \in C^{1,2}([0,T]\times\mathbb{R}^d)$ such that $(t,x)$ is a global minimum (resp. maximum) of $(u_i-\varphi)$, we have,
 \beqs
 && \hspace{-7mm}\mbox{if $t<T$, }
\min\hspace{-1mm}\Big[ \hspace{-1mm}- \Dt{\varphi} -  \Lc^i\varphi   - {f}(i,.,(u_j)_{j\in\Ic},\sigma^\top(i,.) D_x\varphi), \min_{j\in \Ic} {h}_{}(i,j,.,u_i,u_j,\sigma^\top(i,.) D_x \varphi) \Big](t,x)  \geq (\mbox{resp.} \leq) \,  0 ,\\
 && \hspace{-7mm}\mbox{if $t=T$, }
 \min\Big[ u_i-g(i,.) , \; \min_{j\in \Ic}  h_{}(i,j,.,u_i,u_j,\sigma^\top(i,.) D_x \varphi)\Big] (T,x)  \; \geq (\mbox{ resp. } \leq) \;  0\,.
 \enqs
 A locally bounded vector function $(u_i)_{i\in\Ic}$ on $[0,T)\times\R^d$ is called a viscosity solution to \reff{SVI}-\reff{SVITT} if
$u_*$ and $u^*$ are respectively viscosity supersolution and subsolution to \reff{SVI}-\reff{SVITT}.
 \end{Definition}

\begin{Theorem}  \label{probound}
Under {\bf (H0)}-{\bf (H1)}, the function $v$ is a (discontinuous) viscosity solution to \reff{SVI}-\reff{SVITT}.
\end{Theorem}



\noindent {\bf Proof.}  First, following the proof of Lemma 3.3 and Remark 3.2 in Kharroubi et al. (2010), 
standard estimates on the penalized BSDE \reff{BSDEpen}  lead to
 \beqs
\Eb \Big[\sup_{t\in [0,T]} |Y^{e,n}_t|^2\Big]
  \leq
 C\Big( 1 +  \Eb \Big[ |g_{}(I_T^{e},X_T^{e})|^2 + \int_t^T |X_s^{e}|^2 ds+ \sup_{s\in[0,T]}|\tilde v_{I^e_s}(s,X_s^{e})|^2\Big]  \Big),  \quad e\in E\,.
 \enqs
Combining Fatou's lemma with standard estimates on $X$ and linear growth conditions on $g$ and $\tilde v$, see {\bf (H1)}, we get that $\sup_{t\in [0,T]} |v_{i}(t,x)|^2  \leq C (1 +  |x|^2)$ with $C>0$. Thus, $v$ is locally bounded.

We observe that the viscosity property of $v$ in the interior of the domain is based on the same arguments as the one presented in the proof of Theorem 4.1 of Kharroubi et al. (2010). 
 The only difference comes from the more general form of the coefficients $f$ and $h$. This is not a relevant issue here  since they are continuous.
 In order to alleviate the presentation of the paper, we choose to omit it here and  only prove the viscosity property \reff{SVITT} on the maturity boundary.
 
\noindent (i) Let us first consider the supersolution property of $v_*$ to \reff{SVITT}. Let $(i,x_0)\in\Ic\times\mathbb{R}^d$ and $\varphi\in C^{1,2}([0,T]\times\mathbb{R}^d)$ such that $(T,x_0)$ is a null global minimum of $([v_*]_i-\varphi)$. Passing to the limit of the viscosity properties of the penalized BSDE, we
get 
 \beqs
 \min_{j\in \Ic} h_{}(i,j,x,[v_{*}]_i,[v_{*}]_j, \sigma^\top(i,.) D_x \varphi)](T,x_0) &\geq & 0\;.
 \enqs
 Furthermore $v^{n}(T,.)=g$, $n \in \N$, so that the monotonic property of the sequence of continuous functions $(v^n)_{n\in\N}$ gives $v_*(T,.) \ge g$. Therefore $v_*$ is a viscosity supersolution of \reff{SVITT}.
\vspace{2mm}

\noindent (ii) We now turn to the subsolution property of $v^*$. We argue by contradiction and suppose the existence of $(i,x_{0})\in\Ic\times\mathbb{R}^d$ and  $\varphi\in C^{1,2}([0,T]\times\mathbb{R}^d)$ such that
 \beq \label{varphimax}
 0 \; = \; (v_i^*-\varphi)(T,x_{0}) & =  & \max_{[0,T]\times\mathbb{R}^d}(v_i^*-\varphi)\,,
 \enq
and $\min\left[ \varphi-g_{}(i,.) \; , \; \min_{j\in\Ic} h_{}(i,j,.,\varphi,v^*_j, \sigma^\top(i,.) D_x \varphi) \right] (T,x_0)  =:  2\varepsilon \; > \; 0$. The regularity of $v^*$, $\varphi$ and $D_x \varphi$ as well as the monotonic property of $h$ lead to  the existence of an open neighborhood $\Oc$ of $(T,x_{0})$ $\in$ $[0,T]\times\mathbb{R}^d$,
 and $\Upsilon,r>0$ such that for all $(t,x,\eta,\eta')\in \Oc\times(-\Upsilon,\Upsilon)\times B(0,r)$,  we get
 \beq\label{cond0}
 \min\Big[\varphi - \eta - g(i,.) \; , \;\min_{j\in \Ic} h_{}(i,j,.,\varphi-\eta,v^*_j, \sigma^\top(i,.) [D_x \varphi +\eta'])\Big](t,x) & \geq & \eps\;.
 \enq
 We introduce  a sequence $(t_{k},x_{k})_{k}$ valued in $[0,T)\times\R^d$ satisfying $(t_{k},x_{k})\rightarrow(T,x_{0})$ and $v_i(t_{k},x_{k})\rightarrow v_i^*(T,x_{0})$. Let us choose $\delta>0$ such that $[t_{k},T]\times B(x_{k},\delta)\subset \Oc$ for $k$ large enough, and introduce the modified test function $\varphi^{k}$ given by
 \beqs
 \varphi^{k}(t,x) & := & \varphi(t,x)+\left(\zeta\frac{| x-x_{k} |^2}{\delta^2}+C_{k}\phi\left(\frac{x-x_{k}}{\delta}\right)+\sqrt{T-t}\right)\,,
 \enqs
 where $0<\zeta<\Upsilon\wedge\delta r$, $\phi$ is a regular function in $ C^2(\mathbb{R}^d)$ such that $\phi|_{\bar{B}(0,1)}\equiv
0,~\phi|_{\bar{B}(0,1)^{c}}> 0$, $\lim_{|x|\rightarrow\infty}\frac{\phi(x)}{1+|x|}=\infty$, and $C_{k}>0$ is a constant to be determined later. 
 Since $(v^*-\varphi^{k})(t,x) \leq (v^*-\varphi)(t,x)-\zeta\frac{|x-x_{k}|^2}{\delta^2}$ for $(t,x)\in[t_{k},T]\times\R^d$, we deduce from \reff{varphimax} that  $(v^*-\varphi^{k})(t,x)  \leq  -\zeta$, for $(t,x)\in[t_{k},T]\times  \partial B(x_{k},\delta)$. Choosing $C_k$ large enough, the particular form of the function $\phi$ leads to
 \beq  \label{majhdom}
 (v_i^*-\varphi^{k})(t,x) &\leq&  -\frac{\zeta}{2}\;,  \; \quad\mbox{ for  } (t,x)\in  B(x_{k},\delta)^c\times[t_{k},T]\,.
 \enq
 Thanks to the $\sqrt{T-t}$ term in the modified test function $\varphi^{k}$, we deduce that
 \beq
 \left[ - \Dt{\varphi^k} - \mathcal{L}^{i}\varphi^{k} - {f}_{}\Big(i,.,(v^*_j + [\varphi^{k}-\eta - v^*_i]\mathbf{1}_{j=i})_{j\in\Ic},\sigma^\top(i,.) D_x\varphi^{k}\Big)\right](t,x) &\geq&  0 \,, \label{cond1}
 \enq
 for any $(t,x,\eta)\in[t_{k},T)\times B(x_{k},\delta)\times(-\Upsilon+\zeta,\Upsilon)$ and $k$ large enough.
 We now choose $\eta<\Upsilon\wedge\frac{\zeta}{2}\wedge \varepsilon$  and introduce the stopping time
 $
 \theta_{k} :=  \inf\left\{ s\geq t_{k}~;~X_{s}^{e_k} \notin B(x_k,\delta) \;\;\mbox{ or }\;\; I_{s}^{e_k}\neq I_{s-}^{e_k} \right\} \wedge T\,,
 $
 where $e_k:=({t_{k},i,x_{k}})$. Let us finally consider the process $(Y^k,Z^k,U^k,K^k)$ given on $[t_{k},\theta_{k}]$ by
 \beqs
 \left\{
 \begin{array}{rcl}
 Y_{s}^k &:=& \Big[\varphi^{k}(s,X^{e_k}_{s})-\eta\Big]\mathbf{1}_{s\in[t_{k},\theta_{k})} + v_{I^{e_k}_{s}}(\theta_{k},X^{e_k}_{\theta_{k}})\mathbf{1}_{s=\theta_{k}}\,,\quad
 Z^k_{s} \;:=\; \sigma_{}^\top(I^{e_k}_{s-},X^{e_k}_{s}) D_x\varphi^{k}(s,X^{e_k}_{s}),\;\\
 U^k_{s} &:= & \left( \Big[ v^*_j(s,X^{e_k}_{s}) - [\varphi^{k}(s,X^{e_k}_{s})-\eta]\Big] \mathbf{1}_{j\neq {I^{e_k}_{s-}}} \right)_{j\in\Ic}\;,\\
 K^k_{s} & := & - \int_{t_{k}}^{s}\left[\left(\Dt{\varphi^k} + \mathcal{L}^{I^{e_k}_{r}}\varphi^{k}\right)
 + {f}_{}(I^{e_k}_{r},., (v^*_j + [\varphi^{k}-\eta - v^*_{I^{e_k}_{r}}]\mathbf{1}_{j=I^{e_k}_{r}})_{j\in\Ic}, Z^k_r)\right](r,X^{e_k}_r) dr \\
 & & - \; \int_{t_{k}}^s\int_{\Ic}( \varphi^{k}-\eta -v^*_j)(r,X^{e_k}_{r})\mu(dr,dj) + \; \left[\varphi^{k}-\eta-v_{I^{e_k}_{\theta_{k}}}\right](\theta_{k},X^{e_k}_{\theta_{k}})\mathbf{1}_{s=\theta_{k}}.
 \end{array}
 \right.
 \enqs

 One easily checks from \reff{cond0}-\reff{majhdom}-\reff{cond1} that  $(Y^k,Z^k,U^k,K^k)$ is solution to
  \beqs
 Y_{s} =  v_{I_{\theta_k}^{e_k}}(\theta_k, X_{\theta_k}^{e_k})+ \hspace{-1mm}\int_{s}^{\theta_k} \hspace{-2mm} f_{}(I^{e_k}_{r},X^{e_k}_{r},Y_{r}+U_r,Z_{r})dr   -\hspace{-1mm}\int_{s}^{\theta_k}\hspace{-2mm} Z_{r} \cdot dW_{r} - \hspace{-1mm}\int_{s}^{\theta_k}\hspace{-2mm} \int_{\Ic}\hspace{-1mm} U_{r}(j) \mu(dr,dj) + K_{\theta_k}-K_r
 \enqs
 on $[t_k,\theta_k]$, 
 together with the constraint $ h_{}(I_{r-}^{e_k},j,X_{r}^{e_k},Y_{r-}, Y_{r-}+ U_{r}(j),Z_{r})\geq 0$ a.e., $j\in\Ic$. Since $(Y^{e_k},Z^{e_k},U^{e_k},K^{e_k})$ is a minimal solution to this constrained BSDE with jumps, we deduce
 $$
 \varphi^{k}(t_{k},x_{k})-\eta \; = \; \varphi(t_{k},x_{k})+\sqrt{T-t_{k}}-\eta \; \geq \; v_i(t_{k},x_{k})\,, \qquad \mbox{ for all } k \mbox{ large enough. }
 $$
 Letting $k$ go to infinity, this contradicts \reff{varphimax} and concludes the proof. \ep
 
  \vspace{2mm}

\begin{Remark}\label{RemSufficientH1}
 {\rm The main drawback of this representation is the necessity of Assumption {\bf (H1)}. Following similar arguments as in the proof of Proposition 6.3 in Kharroubi et al. (2010), 
 observe that it is satisfied whenever there exists a Lipschitz function $(w_i)_{i\in\Ic}\in [C^{2}(\mathbb{R}^{d} )]^{\Ic}$ supersolution to \reff{SVITT} satisfying a linear growth condition, and there exists a constant $C>0$ such that $\Lc^i w_i + f(.,(w_j)_{j\in\Ic},\sigma_i^\top Dw_i) \leq C $ on $\R^d$, $i\in\Ic$.}
\end{Remark}


\subsection{A comparison argument}\label{subsecComparison}

 In this section, we provide sufficient conditions characterizing the value function $v$ as the unique viscosity solution of \reff{SVI}-\reff{SVITT}. This gives in particular the continuity of $v$, leading to the strong convergence by penalization and the minimality condition, presented in Section \ref{SecBSDE}.
  The proof relies as usual on a comparison argument, which holds under the following additional assumptions.\\[-3mm]
 


\noindent {\bf (H2)} \hspace{1mm} The following holds:\vspace{-2mm}
 \begin{enumerate}[(i)]
\setlength\itemsep{0pt}
\setlength\leftmargin{\parindent}
\setlength\itemindent{-\parindent}
 \item For any $i\in\Ic$, $f(i,.)$ is convex in $((y_j)_{j\in\Ic},z)$ and increasing in $u_i$.
 \item For any $i,j\in\Ic$, $h(i,j,.)$ is concave in $(y_{i},y_{j},z)$ and decreasing in $y_i$. 
\item There exists a nonnegative vector function $(\Lambda_i)_{i\in\Ic}\in[C^2(\mathbb{R}^d)]^\Ic$ and a positive constant $\rho$ such that, for all $i\in\Ic$, $\Lambda_i \geq g_i$,  $\lim_{|x|\rightarrow\infty}\frac{\Lambda_i(x)}{1+|x|}=\infty$ and we have :
 \beqs
 \Lc^i\Lambda_i + f(i,.,(\Lambda_j)_{j\in\Ic},\sigma^\top(i,.) D_x\Lambda_i) \leq \rho\Lambda_i\quad \mbox{ and } \quad
 \min_{j\in \Ic} h_{}(i,j,.,\Lambda_i,\Lambda_j,\sigma^\top(i,.) D_{x}\Lambda_i)>0\;.
 \enqs
  \end{enumerate}

 \vspace{-2mm}
 
\noindent An example where {\bf (H2)} holds is given for the case of optimal switching in  Bouchard (2009).
 \begin{Remark} {\rm  As in Bouchard (2009), 
 (iii) allows us to construct a nice strict supersolution of \reff{SVI}  allowing to control solutions of \reff{SVI}-\reff{SVITT} by convex perturbations. Following the approach of Kharroubi et al. (2010), 
 the general form of $f$ and $h$ forces us to add the extra convexity assumptions (i) and (ii).}
  \end{Remark}

\begin{Theorem} \label{theouni}
Let {\bf (H0)}-{\bf (H1)}-{\bf (H2)} holds. Then, for any $U$ lsc (resp. $V$ usc) viscosity supersolution (resp. subsolution) of  \reff{SVI}-\reff{SVITT} satisfying 
$[|U| + |V|](t,x) \leq C(1+|x|)$ on $[0,T]\times\R^d$, we have $U_i\geq V_i$ on $[0, T ]\times\R^d$, $i\in\Ic$.
In particular, $v$ is continuous and it is the unique viscosity solution of \reff{SVI}-\reff{SVITT} satisfying a linear growth condition.
\end{Theorem}

 We omit the proof of this comparison theorem which is a natural extension of Theorem 4.1 in Kharroubi et al. (2010). 
  Following the arguments of the proof of Proposition 3.3 in Peng and Xu (2007), 
  $v$ can still be  interpreted as the minimal viscosity solution of \reff{SVI}-\reff{SVITT} in the class of functions with linear growth, whenever a comparison theorem for the IPDE \reff{IPDE} holds.



\section{Numerical issues}\label{SectionNumerics}

\setcounter{equation}{0} \setcounter{Assumption}{0}
\setcounter{Theorem}{0} \setcounter{Proposition}{0}
\setcounter{Corollary}{0} \setcounter{Lemma}{0}
\setcounter{Definition}{0} \setcounter{Remark}{0}

 The numerical resolution of systems of variational inequalities of the form \reff{SVI}-\reff{SVIT} usually relies on the use of iterated free boundary. 
 We first solve the system without boundary condition and consider recursively the system constrained by the boundary condition coming from the previous iteration. 
 In a switching problem, we constrain the solution associated to $n+1$ possible switches by the obstacle built from the solution where only $n$ switches are allowed. Such a numerical approach is computationally demanding.  
 We present here a natural convergent algorithm based on the approximation of the solution to the corresponding constrained BSDE with jumps \reff{BSDEgen}-\reff{cons}.
 We combine a penalization procedure with the discrete-time scheme studied by Bouchard and Elie (2008) 
 and the statistical estimation projection presented in Gobet et al. (2006). 
 Thanks to the previous Feynman-Kac representation, this gives rise to a convergent probabilistic algorithm solving coupled systems of variational inequalities.
\vspace{2mm}

We fix an initial condition $e\in E$ and omit it in the expressions for ease of presentation. Suppose that {\bf (H0)}-{\bf (H1)}-{\bf (H2)} holds. The algorithm is divided in three steps.

 \vspace{1mm}

\noindent {\bf Step 1. Approximation by penalization.} We first approach the constrained BSDE with jumps \reff{BSDEgen}-\reff{cons} by its penalized version \reff{BSDEpen} characterized by a driver $f^n:=f-n [h]^-$ as in Section \ref{SubSecPenalization}. We deduce from Proposition \ref{PropConvPenError} that the penalization error converges to $0$ as $n$ goes to infinity, see \reff{ConvPenError2}.


\vspace{1mm}

\vspace{1mm}

\noindent {\bf Step 2. Time discretization.}  Observe that the pure
jump process $I$ can be simulated perfectly and denote by $(\tau_l)_{l}$
its jump times on $[0,T]$. We introduce the Euler time scheme
approximation $X^h$ of the forward process $X$ defined on the
concatenation $(s_l)_l$ of the regular time grid $\{t_k:= kh , \;
k=1,\ldots, T/h  \}$ with the jumps $(\tau_l)$ of $I$: 
 \beqs
 X^{h}_0= X_0 &\mbox{ and }&
 X^{h}_{s_{l+1}} := X^{h}_{s_l} + b_{}(I_{s_l},X^{h}_{s_l}) (s_{l+1}-s_{l}) + \sigma_{}(I_{s_l},X^{h}_{s_l}) [W_{s_{l+1}}-W_{s_{l}}] .
 \enqs
 We deduce an approximation $Y^{n,h}_T$ of $Y^n_T$ at maturity given by $g_{I_T}(X^h_T)$.
 The penalized BSDE \reff{BSDEpen} can now be discretized by an extension of the scheme exposed in Bouchard and Elie (2008)  
 An approximation of $Y^n$ at time 0 is computed recursively following the backward scheme for $k=T/h - 1,\cdots,0$ :
   \beq
   \left\{
   \begin{array}{lll}
   Z^{n,h}_{t_k} &:= &{1\over h} \Eb_{t_k}\Big[ Y^{n,h}_{t_{k+1}}(W_{t_{k+1}}-W_{t_k})\Big] \label{Znh}\\
   U^{n,h}_{t_k}(i) &:=& {1\over h}\Eb_{t_k}\Big[Y^{n,h}_{t_{k+1}}\frac{\tilde\mu((t_k,t_{k+1}]\times \{i\})}{\lambda(i)}\Big], \quad i\in\Ic \label{Unh}\\
   Y^{n,h}_{t_k} &:=& \Eb_{t_k}\left[ Y^{n,h}_{t_{k+1}} + \int_{t_{k}}^{t_{k+1}} f^n_{}(I_{s},X^{h}_{t_k}, Y^{n,h}_{t_{k+1}}, Z^{n,h}_{t_k},U^{n,h}_{t_k}) ds \right]\label{Ynh}
   \end{array}
   \right.
   \enq
 where $\Eb_{t_k}$ denotes the conditional expectation with respect to $\Gc_{t_k}$.  Following the arguments of Section 2.5 in Bouchard and Elie (2008) 
 and identifying $(Y^{n,h},Z^{n,h},U^{n,h})$ as a process constant on each interval $(t_k,t_{k+1}]$, we verify the convergence of this discrete-time approximation~:
 \beq\label{ConvDiscError}
  \|Y^{n}-Y^{n,h}\|_{_{\bf \Sc^2}}   +  \|Z^{n}-Z^{n,h}\|_{_{\bf L^2_{W}}}  +  \|U^{n}-U^{n,h}\|_{_{\bf L^2_{\tilde\mu}}}   
  & \mathop{{\longrightarrow}}\limits_{h\rightarrow\infty} &  0,\quad n\in\N.
 \enq

 \noindent  {\bf Step 3. Approximation of the conditional expectations.}  
 The last step consists in estimating the conditional expectation operators $\Eb_{t_k}$ arising in \reff{Znh}. 
 We adopt here the approach of Longstaff-Schwarz generalized in Gobet et al. (2006) 
 relying on least square regressions. 

 Fix $N\in\N$ and simulate $N$ independent copies of the Brownian increments $(W^j_{t_{k+1}}-W^j_{t_k})_{0\le k \le T/h}$ and the poisson measure $(\tilde\mu^j ((t_k,t_{k+1}]\times\Ic)_{0\le k \le T/h}$. For each simulation $j\le N$, define $I_j^{N}$ and $X_j^{h,N}$ as the trajectories of $I$ and $X^h$.
 By induction, one can easily verify the Markov property of the process $(Y^{n,h},Z^{n,h},U^{n,h})$ defined in \reff{Znh}:
 $$
 Y^{n,h}_{t_i} = c^{n,h}_k(I_{t_i},X^{h}_{t_i}),\quad
 Z^{n,h}_{t_i} = a^{n,h}_k(I_{t_i},X^{h}_{t_i}), \quad
 U^{n,h}_{t_i} = b^{n,h}_k(I_{t_i},X^{h}_{t_i}),
 $$
 for some deterministic functions $(a^{n,h}_k, b^{n,h}_k, c^{n,h}_k)_{k\le n}$.
 The idea is to approximate these functions using Ordinary Least Square (OLS) estimators.
 Given $L\in\N$, we introduce a collection of basis functions $(a^L_l, b^L_l, c^L_l)_{1\le l \le L}$ of $\R\times\R^d\times\R^d$.
 For each trajectory $j\le N$, define the associated terminal value given by $Y^{n,h,L,N}_{j,t_n} := g_{I^N_{k,t_n}}(X^{h,N}_{j,t_n})$.
 Now we define recursively  $(Z^{n,h,L,N}_{j,t_k}, U^{n,h,L,N}_{j,t_k})$, backward in time for $k=T/h-1,\cdots,0$, by computing the OLS approximations  as follows:
 \beqs 
 (\hat\alpha_1,\cdots,\hat\alpha_L) :=
 {\rm arg}\min_{\alpha_1,\cdots,\alpha_L} {1\over N}\sum_{j=1}^N \Big|{1\over h} Y^{n,h,L,N}_{j,t_{k+1}}[W^j_{t_{k+1}}-W^j_{t_k}] - \sum_{l=1}^L \alpha_l a^L_l(I^{N}_{j,t_k}, X^{h,N}_{j,t_k})\Big|^2,\\
 (\hat\beta_,\cdots,\hat\beta_L)(i) :=
 {\rm arg}\min_{\beta_1,\cdots,\beta_L} {1\over N} \sum_{j=1}^N \Big|{1\over h} Y^{n,h,L,N}_{t_{k+1}}\frac{\tilde\mu^j((t_k,t_{k+1}]\times {\{i\}})}{\lambda(i)} - \sum_{l=1}^L \beta_l b^L_l(I^{N}_{j,t_k}, X^{h,N}_{j,t_k})\Big|^2\,,
 \enqs
 for $i\in\Ic$, leading to the approximation
 \beqs
 &&
 Z^{n,h,L,N}_{j,t_k}:= \sum_{l=1}^L \hat\alpha_l a^L_l(I^{N}_{j,t_k}, X^{h,N}_{j,t_k})
 \;\;\mbox{ and }\;\;
 U^{n,h,L,N}_{j,t_k}(i):= \sum_{l=1}^L \hat\beta_l(j) b^L_l(I^{N}_{j,t_k}, X^{h,N}_{j,t_k})\,,\;\; i\in\Ic.
 \enqs
It remains to introduce $(\hat\gamma_1,\cdots,\hat\gamma_L)$ the minimizer  of the mean square error
 \beqs
{1\over N}\sum_{j=1}^N \Big|Y^{n,h,L,N}_{j,t_{k+1}}    + \int_{t_{k}}^{t_{k+1}} \hspace{-2mm} f^n_{}(I^{N}_{j,s},X^{h,N}_{j,t_k}, Y^{n,h,L,N}_{t_{k+1}}, Z^{n,h,L,N}_{t_k},U^{n,h,L,N}_{t_k}) ds    - \sum_{l=1}^L \gamma_l c^L_l(I^{N}_{j,t_k}, X^{h,N}_{j,t_k})\Big|^2
 \enqs
 in order to deduce the OLS approximation $Y^{n,h,L,N}_{j,t_k} := \sum_{l=1}^L \hat\gamma_l c^L_l(I^{N}_{j,t_k}, X^{h,N}_{j,t_k})$.
 
 We refer to Gobet et al. (2006) 
 for the control of the statistical error due to the approximation of the conditional expectation operators by OLS projections, and, by extension, 
 \beq\label{ConvStatError}
  \|Y^{n,h}-Y^{n,hL,N}\|_{_{\bf \Sc^2}}   +  \|Z^{n,h}-Z^{n,h,L,N}\|_{_{\bf L^2_{W}}}  +  \|U^{n,h}-U^{n,h,L,N}\|_{_{\bf L^2_{\tilde\mu}}}
   \mathop{{\longrightarrow}}\limits_{N,L\rightarrow\infty}     0\,, \; n\in\N\,,\, h>0.\;\;
 \enq
The convergence of the algorithm follows from \reff{ConvPenError2}, \reff{ConvDiscError} and \reff{ConvStatError}. The derivation of a convergence rate requires precisions on the influence of $n$ on the discretization and statistical errors, as well as a control of the penalization error. This challenging point is left to further research.

\end{document}